\documentclass[12pt,reqno]{amsart}
\usepackage{pifont}
\usepackage{mathrsfs}
\usepackage{geometry}
\usepackage{titletoc}
\usepackage{stix2}

\usepackage{amsmath}
\usepackage{amssymb} 
\usepackage{enumitem} 
\usepackage{mathtools} 
\usepackage[table]{xcolor} 
\usepackage[all]{xy} 
\usepackage{tikz} 
\usepackage{tikz-cd}
\usepackage{indentfirst} 
\usepackage{babel} 
\usepackage{setspace} 

\usepackage[colorlinks,linkcolor=red,anchorcolor=green,citecolor=blue]{hyperref} 
\hypersetup{linktocpage = true} 

\usepackage{rotating} 

\usepackage{ytableau} 
\usepackage{longtable} 
\newcolumntype{M}[1]{>{\centering\arraybackslash}m{#1}} 

\usepackage{fancyhdr}

\newcommand\cO{{\mathcal O}}

\newcommand\cX{{\mathcal X}}

\newcommand\sR{{\mathcal R}}

\newcommand{\id}{{\rm id}}

\DeclareMathOperator*{\Sym}{Sym}

\newcommand\bp{{\bar\partial}}

\theoremstyle{plain}
\newtheorem{thm}{Theorem}[section]
\newtheorem{lemma}[thm]{Lemma}
\newtheorem{prop}[thm]{Proposition}
\newtheorem{cor}[thm]{Corollary}
\newtheorem{defn}[thm]{Definition}

\theoremstyle{definition}
\newtheorem{example}[thm]{Example}
\newtheorem{remark}[thm]{Remark}

\newcommand{\btheorem}{\begin{thm}}
    \newcommand{\etheorem}{\end{thm}}
\newcommand{\bproposition}{\begin{prop}}
    \newcommand{\eproposition}{\end{prop}}
\newcommand{\bdefinition}{\begin{defn}}
    \newcommand{\edefinition}{\end{defn}}
\newcommand{\bcorollary}{\begin{cor}}
    \newcommand{\ecorollary}{\end{cor}}
\newcommand{\bproof}{\begin{proof}}
    \newcommand{\eproof}{\end{proof}}
\newcommand{\bremark}{\begin{remark}}
    \newcommand{\eremark}{\end{remark}}
\newcommand{\eexample}{\end{example}}
\newcommand{\bexample}{\begin{example}}

\newcommand{\elemma}{\end{lemma}}
\newcommand{\blemma}{\begin{lemma}}

\newcommand{\la}{\langle}
\newcommand{\ra}{\rangle}
\newcommand{\sq}{\sqrt{-1}}

\newcommand{\p}{\partial}

\renewcommand{\bar}{\overline}

\newcommand{\beq}{\begin{equation}}
\newcommand{\eeq}{\end{equation}}
\newcommand{\ee}{\end{eqnarray*}}
\newcommand{\be}{\begin{eqnarray*}}

\newcommand{\bd}{\begin{enumerate}}
    \newcommand{\ed}{\end{enumerate}}

\renewcommand{\hat}{\widehat}
\renewcommand{\tilde}{\widetilde}

\newcommand{\qtq}[1]{\quad\mbox{#1}\quad}
\renewcommand{\bp}{\bar{\partial}}
\newcommand{\Om}{\Omega}

\renewcommand{\>}{\rightarrow}

\newcommand{\C}{{\mathbb C}}

\renewcommand{\P}{{\mathbb P}}

\newcommand{\R}{{\mathbb R}}

\newcommand{\LL}{\left\langle}
\newcommand{\RL}{\right\rangle}

\setlist[itemize]{leftmargin=*}
\setlist[enumerate]{leftmargin=*}

\numberwithin{equation}{section} 

\makeatletter

\usepackage{fancyhdr}
\newcommand{\om}{\omega}
\renewcommand{\le}{\leqslant}
\renewcommand{\leq}{\leqslant}

\renewcommand{\geq}{\geqslant}
\newcommand{\vone}{\vskip 1\baselineskip}
\newcommand{\Ric}{\operatorname{Ric}}

\title{Laplacian comparison theorems on complete K\"ahler manifolds and applications}

\author{Jiaxuan Fan}
\author{Zhiyao Xiong}
\author{Xiaokui Yang}

\address{Jiaxuan Fan, Qiuzhen College, Tsinghua University, Beijing, 100084, China}
\email{fanjiaxu21@mails.tsinghua.edu.cn}

\address{Zhiyao Xiong, Department of Mathematics, Tsinghua University, Beijing, 100084, China}
\email{xiongzy22@mails.tsinghua.edu.cn}

\address{Xiaokui Yang, Department of Mathematics and Yau Mathematical Sciences Center, Tsinghua University, Beijing, 100084, China}
\email{xkyang@mail.tsinghua.edu.cn}

\begin{document}
\begin{abstract}
In this paper, we establish new Laplacian comparison theorems and rigidity theorems for complete K\"ahler manifolds under new curvature notions that interpolate between Ricci curvature and holomorphic bisectional curvature.
\end{abstract}
\maketitle
    \setcounter{tocdepth}{1}

{\small{    \begin{spacing}{1.1} \tableofcontents %
            \dottedcontents{section}[1.8cm]{}{3em}{5pt} %
\end{spacing} }}

\section{Introduction}
Comparison theorems serve as essential analytical tools in Riemannian geometry, unveiling profound connections between curvature bounds and the geometric-topological structure of manifolds.
Let $(M,g)$ be a complete $n$-dimensional Riemannian manifold with Ricci curvature $\mathrm{Ric}(g)\geq (n-1)kg$.
The foundational Laplacian comparison theorem (see, e.g., \cite{Pet16, Lee18}) asserts that for any fixed point $p \in M$, the distance function $r(x) = d(p,x)$ satisfies
\begin{equation}
\Delta r(x)
\leq (n-1)\frac{\operatorname{sn}_k'(r(x))}{\operatorname{sn}_k(r(x))}
\label{intro Laplacian}
\end{equation}
on $M \setminus \operatorname{cut}(p) \cup \{p\}$. Notably, the identity in \eqref{intro Laplacian}  holds if and only if the universal cover of $(M,g)$ is isometric to a space form.
The Laplacian comparison theorem serves as a fundamental tool for deriving numerous comparison results in Riemannian geometry. For instance, it implies Myers' diameter theorem, which states that if $\mathrm{Ric}(g) \geq (n-1)g$, then $\mathrm{diam}(M,g) \leq \mathrm{diam}(\mathbb{S}^n, g_{\mathrm{can}}) = \pi$. Building upon the Laplacian comparison, Bishop and Gromov (e.g. \cite{BC64}, \cite{CE08}) established the volume comparison theorem, a powerful tool in global geometry. This result, in turn, plays a crucial role in proving Cheng's diameter rigidity theorem \cite{Che75}, which asserts that if $\mathrm{Ric}(g) \geq (n-1)g$ and $\mathrm{diam}(M,g) = \pi$, then $(M,g)$ is isometric to the round sphere $(\mathbb{S}^n, g_{\mathrm{can}})$.
For more details along this comprehensive topic, we refer to \cite{CC97},  \cite{Zhu97} and \cite{Wei07} and the references therein.\\

In K\"ahler geometry, the model spaces for comparison are the simply connected K\"ahler manifolds $M_c$ with constant holomorphic bisectional curvature $\mathrm{HBSC}\equiv c$.
Li and Wang  established in \cite{LW05} the Laplacian comparison theorem that if $\mathrm{HBSC}\geq c$, then for any fixed point $p \in M$, the distance function $r(x) = d(p,x)$ satisfies
\begin{equation}
\Delta r\leq 2(n-1)\frac{\mathrm{sn}_{c/2}'(r)}{\mathrm{sn}_{c/2}(r)}+\frac{\mathrm{sn}_{2c}'(r)}{\mathrm{sn}_{2c}(r)}
\label{intro LW Laplacian comparison}
\end{equation}
on $M \setminus \operatorname{cut}(p) \cup \{p\}$. This estimate also enjoys rigidity: the identity in \eqref{intro LW Laplacian comparison} holds if and only if the universal cover of $(M,\om_g)$ is isometrically biholomorphic to $M_c$.
Moreover, they proved that if $\mathrm{HBSC}\geq 1$, then $\mathrm{diam}(M,\om_g)\leq \mathrm{diam}(\C\P^n, \omega_{\mathrm{FS}})$ and $\mathrm{Vol}(M,\om_g)\leq \mathrm{Vol}(\C\P^n, \omega_{\mathrm{FS}})$, with the identity in the volume comparison holding if and only if $(M, \omega_g)$ is isometrically biholomorphic to $(\C\P^n, \omega_{\mathrm{FS}})$. More recently, Datar and Seshadri \cite{DS23} established the diameter rigidity theorem, which states that if $\mathrm{HBSC}\geq 1$ and $\mathrm{diam}(M,g)=\mathrm{diam}(\C\P^n, \omega_{\mathrm{FS}})$, then $(M,\omega_g)$ is isometrically biholomorphic to $(\C\P^n,\omega_{\mathrm{FS}})$. This is achieved by using Siu-Yau's solution to the Frankel conjecture \cite{SY80} and an interesting monotonicity formula for Lelong numbers on $\C\P^n$ (\cite{Lot21}).
To the best of our knowledge, this approach has no counterpart in classical Riemannian geometry.\\

 Utilizing entirely different techniques in algebraic geometry (e.g. \cite{Fuj18}), Zhang \cite{Zha22} remarkably established volume comparison and rigidity theorems under Ricci lower bounds:

\btheorem\label{Zhang} Let $(M^n,\omega_g)$ be a complete K\"ahler manifold with $\mathrm{Ric}(\omega_g) \geq (n+1)\omega_g$. Then
 \beq \mathrm{Vol}(M,\omega_g)\leq \mathrm{Vol}(\C\P^n, \omega_{\mathrm{FS}}).\label{ZKW}\eeq
Moreover, the identity in \eqref{ZKW} holds if and only if $(M,\omega_g)$ is isometrically biholomorphic to  $(\C\P^n, \omega_{\mathrm{FS}})$.
\etheorem

\noindent
One might naturally consider comparing the diameter $d$ and the Laplacian $\Delta_M r$ of such manifolds with those of the model space $(\C\P^n, \omega_{\mathrm{FS}})$.  However, unlike the case where $\mathrm{HBSC}\geq 1$,  neither of these comparisons holds when the complex dimension $n\geq 2$. This is clearly illustrated by the example  $\left(\C\P^{1},\frac{2}{3}\om_{\mathrm{FS}}\right)\times \left(\C\P^{1},\frac{2}{3}\om_{\mathrm{FS}}\right).$  Furthermore, we demonstrate that the model Laplacian comparison \eqref{intro LW Laplacian comparison}  may fail even locally (see Example \ref{examplen=2}). For more related comparison theorems, we refer to \cite{LW05,Mun09, Liu11, TY12, Liu14, LY18, NZ18, Zhu22, CLZ24+, XY24+, Yang25+} and the references therein.\\

In this paper, we establish new Laplacian comparison theorems and applications for complete K\"ahler manifolds under new curvature notions that interpolate between Ricci curvature and holomorphic bisectional curvature. Let $(M,\omega_g)$ be a K\"ahler manifold, and let $R$ denote its Chern curvature.
One can define the \emph{symmetrized curvature operator}  $ \mathcal {R}: \Gamma(M,\mathrm{Sym}^2 T^{1,0}M) \> \Gamma(M,\mathrm{Sym}^2 T^{1,0}M )$ by the relation \beq  g\left( \mathcal{R} \left(a\right),  b\right) = R_{i\bar j k \bar\ell} a^{ik}\bar b{}^{j\ell} \eeq where  $a=\sum a^{ik}\frac{\p}{\p z^i} \otimes \frac{\p}{\p z^k}$ and $b=\sum b^{j\ell}\frac{\p}{\p z^j} \otimes \frac{\p}{\p z^\ell}$ are in  $\Gamma(M,\mathrm{Sym}^2 T^{1,0}M)$ (see \cite{CV60}, \cite{BNPSW25} and \cite{WY25+}) . This curvature notion appears naturally in various Bochner-Kodaira formulas. For instance, Wang and the third named author established in \cite[Theorem~1.2]{WY25+} new Bochner-Kodaira  formulas with quadratic curvature terms on compact K\"ahler manifolds:  for any $\eta\in \Om^{p,q}(M)$,
\beq   \la\Delta_{\bar \p} \eta,\eta\ra =\la \Delta_{\bp_F} \eta,\eta\ra +\frac{1}{4}\left\langle \left(\mathcal {R} \otimes \mathrm{Id}_{\Lambda^{p+1,q-1}T^*M}\right)(\mathbb T_\eta),\mathbb T_\eta \right\rangle.
\eeq 
This linearized curvature term yields new vanishing theorems and provides estimates for Hodge numbers under exceptionally weak curvature conditions. Please  refer to \cite{WY25+} for further discussion on applications.  We say that  $(M,\omega_g)$ has \emph{positive symmetrized curvature operator} $\mathcal R$
if it is positive definite as a Hermitian bilinear form.
A straightforward computation shows that the symmetrized curvature operator of $(\C\P^n, \omega_{\mathrm{FS}})$ is $\mathcal R=2\cdot \mathrm{id}$ which is  positive definite. On the other hand,
if  $\mathcal {R}$ is a positive operator, then $(M,g)$ has positive holomorphic bisectional curvature. Furthermore,  when $M$ is compact,  it follows from Siu-Yau's solution to the Frankel conjecture (\cite{SY80,Mori1979}) that  $M$ is biholomorphic to $\C\P^n$. The following weaker notion on $k$-positivity is natural:

\bdefinition Let $A$ be a Hermitian $ n \times n $ matrix and $ \lambda_1 \leq \cdots \leq \lambda_n $ be eigenvalues of $ A $.   It  is said to be $ k $-positive if
\beq \lambda_1 + \cdots + \lambda_k > 0. \eeq
The symmetrized curvature operator $ \mathcal{R}: \Gamma(M,\Sym^2 T^{1,0}M) \> \Gamma(M,\Sym^2 T^{1,0}M) $ is called $ k $-positive if $ \mathcal{R} $ is $k $-positive at every point  of $ M $. One can define $k$-semi-positivity, $k$-negativity and $k$-semi-negativity in similar ways.
\edefinition
\noindent  By using a combinatorial computation, one can show that  if $(M,\omega_g)$ is a compact K\"ahler manifold of complex dimension $n$ and $\mathcal R$ is $k$-positive with $k\leq (n+1)/2$, then it has positive Ricci curvature (Corollary \ref{symmetry to ric}). Moreover, if $(M,\omega_g)$ is the hyperquadric in $\C\P^{n+1}$ with the induced metric, then  $ \mathcal{R} $ has eigenvalues   (e.g. \cite{CV60}, \cite{BNPSW25} and \cite{WY25+})  $$ \lambda_1 = 2-n  \qtq{and} \lambda_2 = \cdots = \lambda_{N} = 2 $$  where $ N = \frac{n(n+1)}{2} $. In particular, $\mathcal R$  is $\left(\left\lfloor \frac{n}{2}\right\rfloor+1\right) $-positive.\\

The first main result of this paper is  the global Laplacian
comparison theorem for the symmetrized curvature operator $
\mathcal{R} $ with \emph{negative} lower bounds:

\btheorem\label{thm symmetric curvature operator comparison}
Let $(M,\om_g)$ be a complete K\"ahler manifold of complex dimension $n$ and $r$ be  the distance function from a given point $p\in M$.
If $\mathcal{R}-2c\cdot\mathrm{id} $ is $k$-semipositive for some $c<0$ and $k\leq (n+1)/2$, then
\beq
\Delta r\leq2(n-1)\frac{\mathrm{sn}_{c/2}'(r)}{\mathrm{sn}_{c/2}(r)}+\frac{\mathrm{sn}_{2c}'(r)}{\mathrm{sn}_{2c}(r)}
\eeq
on $M\setminus \mathrm{cut}(p)\cup \{p\}$.
Moreover, if the identity in \eqref{Laplacian comparison 2} holds on $M\setminus\mathrm{cut}(p)\cup \{p\}$,
then the universal cover of $(M,\omega_g)$ is isometrically biholomorphic to the hyperbolic space $M_c$.
\etheorem

\noindent  We also obtain a local Laplacian comparison theorem for
the symmetrized curvature operator $ \mathcal{R} $ with
\emph{positive} lower bounds:

\btheorem \label{symmetric laplacian positive} Let $(M,\om_g)$ be a
complete K\"ahler manifold of complex dimension $n$ and $r$ be  the
distance function from a given point $p\in M$. If $\mathcal{R}-2c\cdot\mathrm{id} $
is $k$-semi-positive for some $c>0$ and $k< (n+1)/2$, then \beq
\Delta
r\leq2(n-1)\frac{\mathrm{sn}_{c/2}'(r)}{\mathrm{sn}_{c/2}(r)}+\frac{\mathrm{sn}_{2c}'(r)}{\mathrm{sn}_{2c}(r)}
\label{Laplacian comparison 0} \eeq on some metric ball
$B(p,C(k,n)/\sqrt{c})\setminus \mathrm{cut}(p)\cup \{p\}$. Moreover
if the identity \eqref{Laplacian comparison 0} holds for all such
$x$, then $\mathrm{HBSC}\equiv c$ on $B(p,C(k,n)/\sqrt{c})$.
\etheorem

\noindent Theorems \ref{thm symmetric curvature operator comparison} and \ref{symmetric laplacian positive} are established through the application of a new index theorem (Theorem~\ref{thm new index}) in combination with a combinatorial curvature synthesis technique. On the other
hand, it is clear that Theorem \ref{thm symmetric curvature operator
comparison} and Theorem \ref{symmetric laplacian positive} can give
Bishop-Gromov type \emph{local} volume comparison theorem. We also
establish a diameter comparison theorem when $c>0$ and $k\leq n$:
\btheorem\label{thm diam} Let $(M, \omega_g)$ be a complete K\"ahler
manifold of complex dimension $n$. If $\sR-2c\cdot\id$ is
$k$-semi-positive for some $c>0$ and $k\leq n$, then $M$ is compact
and \beq \mathrm{diam}(M,\om_g)\leq
\frac{\pi}{\sqrt{\nu}}\qtq{where}\nu=\frac{2kc}{4k-3}.\label{diam
estimate} \eeq Moreover, the fundamental group $\pi_1(M)$ of $M$ is
trivial. \etheorem

\noindent  Recall that the symmetrized curvature operator of the hyperquadric in $\mathbb{C}\mathbb{P}^{n+1}$ is $\left(\left\lfloor \frac{n}{2}\right\rfloor+1\right)$-positive. One can observe that the $k$-positivity condition ($k \leq n$) in Theorem \ref{thm diam} is rather weak. It remains entirely unclear whether such manifolds are Fano or not. Even when $k \leq (n+1)/2$, one has $\Ric(\omega_g)\geq c(n+1)\omega_g$, which can imply a Myers-type diameter estimate using the underlying Riemannian metric. Crucially, however, the estimate in \eqref{diam estimate} is sharper than this Riemannian comparison, as further demonstrated in Section \ref{example}. The key ingredient in the proof of Theorem
\ref{thm diam} is a weighted Hessian estimate, which relies on the
derived mixed curvature estimate. As is standard in Riemannian
geometry, this estimate implies a diameter bound, which in turn
implies the finiteness of $\pi_1(M)$ and $b_1(M) = 0$.  By using
vanishing theorems derived in \cite{WY25+} and the Riemann-Roch
theorem, we conclude that $M$ is simply connected.

\vone
Furthermore,  using similar ideas as in the proofs of Theorem \ref{thm symmetric curvature operator comparison} and Theorem \ref{symmetric laplacian positive}, we obtain the following result under the conditions holomorphic sectional curvature $\mathrm{HSC} \geq c$ and Ricci curvature $\Ric \geq c(n+1)\omega_g$:

\btheorem \label{thm Laplacian Ric and HSC}
Let $(M, \omega_g)$ be a complete K\"ahler manifold of complex dimension $n$.
If $\mathrm{HSC}\geq 2c$ and $\mathrm{Ric}\geq c(n+1)\om_g$ for some $c>0$, then the distance function $r$ from a  given point $p\in M$  satisfies
\beq
\Delta r\leq 2(n-1)\frac{\mathrm{sn}_{c/2}'(r)}{\mathrm{sn}_{c/2}(r)}+\frac{\mathrm{sn}_{2c}'(r)}{\mathrm{sn}_{2c}(r)} \label{intro Laplacian comparison 1}
\eeq
on $M\setminus\mathrm{cut}(p)\cup \{p\}$. Moreover, if the identity in \eqref{intro Laplacian comparison 1} holds on $M\setminus\mathrm{cut}(p)\cup \{p\}$,
then $(M,\omega_g)$ is isometrically biholomorphic to the projective space $M_c$.
\etheorem

\noindent Using the Laplacian Comparison Theorem (Theorem \ref{thm Laplacian Ric and HSC}), we can derive the following local Bishop-Gromov type volume comparison theorem:

\btheorem\label{thm volume Ric and HSC}
Let $(M, \omega_g)$ be an $n$-dimensional  complete K\"ahler manifold with $\mathrm{HSC}\geq 2c$ and $\mathrm{Ric}\geq c(n+1)\om_g$ for some $c>0$.
Suppose that $B(p,\delta)\subset M$ is the metric ball centered at $p\in M$ with radius $\delta$, and ${B}(\tilde{p},\delta)$ is a corresponding metric ball in $M_{c}$.
Then the volume ratio
\beq
\frac{\mathrm{Vol}(B(p,\delta))}{\mathrm{Vol}(B(\tilde{p},\delta))}
\eeq
is non-increasing in $\delta$, and in particular
\beq \mathrm{Vol}(B(p,\delta))\leq  \mathrm{Vol}(B(\tilde p,\delta)).\label{intro volume comparison}\eeq
Moreover, if the identity in \eqref{intro volume comparison} holds for some $\delta>0$, then $\omega_g$ has constant holomorphic bisectional curvature $c$ on $B(p,\delta)$.
\etheorem

\noindent As an application of Theorem \ref{thm volume Ric and HSC}, we provide a \emph{differential geometric proof} of the following volume comparison theorem, which is a special case of Theorem \ref{Zhang}.
\bcorollary \label{cor volume Ric and HSC}
Let $(M,\om_g)$ be a complete K\"ahler manifold of complex dimension $n$.
If $\mathrm{HSC}\geq 2c$ and $\mathrm{Ric}\geq c(n+1)\om_g$ for some $c>0$, then $M$ is compact and  \beq \mathrm{Vol}(M,\omega_g)\leq \mathrm{Vol}(M_c).\label{BGG}\eeq
Moreover, the identity in \eqref{BGG} holds if and only if $(M,\omega_g)$ is isometrically biholomorphic to the projective space $M_c$.
\ecorollary

\noindent\textbf{Acknowledgements}.  The third named author wishes to thank Bing-Long Chen for inspired discussions.

\vskip 2\baselineskip

\section{Index forms on K\"ahler manifolds}
In this section, we introduce a new index form for establishing Hessian or Laplacian comparison theorems on K\"ahler manifolds.
Let $(M,\om_g)$ be a K\"ahler manifold.
For a unit-speed geodesic $\gamma:[0,\ell]\to M$,  we set \beq E_\gamma=\frac{\gamma'-\sq J\gamma'}{\sqrt{2}} \in \Gamma\left([0,\ell], \gamma^*T^{1,0}M\right).\eeq
For any $V,W\in \Gamma\left([0,\ell],\gamma^*T^{1,0}M\right)$,  there is  an index form $\cX_\gamma$ of $\gamma$ given by
\beq
    \cX_\gamma\left(V,\bar W\right)=\int_0^\ell \LL V',W'\RL-\frac{1}{2}R\left(E_\gamma,\bar{E_\gamma},V,\bar W\right)dt,\label{new index form X}
\eeq
where $V'=\hat\nabla_{\frac{d}{dt}}V$ and $\hat\nabla$ is the complexification of the pullback connection on $\gamma^*TM$.
We establish the following version of the Hessian comparison theorem:
\btheorem \label{thm new index}
    Let $(M,\om_g)$ be a complete K\"ahler manifold, and $r$ be the distance function from a fixed point $p\in M$.
    For a given point $x\in M\setminus \mathrm{cut}(p)\cup \{p\}$ and $X\in T_x^{1,0}M$,  if  $\gamma:[0,\ell]\to M$ is the unit-speed minimal geodesic connecting $p$ and $x$, then
    \beq
        \left(\p\bp r\right)\left(X,\bar X\right)\le \cX_\gamma\left(V,\bar V\right)-\frac{1}{2}\int_0^\ell \left|\LL V',E_\gamma\RL\right|^2dt,\label{new index form estimate}
    \eeq
    where $V\in \Gamma([0,\ell],\gamma^*T^{1,0}M)$ is a  vector field  satisfying $V(0)=0$ and $V(\ell)=X$.  Moreover, the identity in \eqref{new index form estimate} holds if and only if $V^\perp=V-\LL V,\gamma'\RL\gamma'$ is a  Jacobi field.
\etheorem

\begin{proof}
    Since $V\in \Gamma([0,\ell],\gamma^*T^{1,0}M)$, we  assume that $V=\frac{1}{\sqrt{2}}\left(\tilde V-\sq J\tilde V\right)$ where $\tilde V$ is a real vector field along $\gamma$.
    If we set $v=\tilde V(\ell)$, then one has
    \beq
        X=V(\ell)=\frac{1}{\sqrt{2}}\left(\tilde V(\ell)-\sq J\tilde V(\ell)\right)=\frac{1}{\sqrt{2}}\left(v-\sq Jv\right).
    \eeq
A simple calculation shows that
    \beq
        2\left(\p\bp r\right)\left(X,\bar X\right)=  \left(\mathrm{Hess}\,r\right)(v,v)+ \left(\mathrm{Hess}\,r\right)(Jv,J v).
    \eeq
    If we write $a=\LL v,\gamma'(\ell)\RL$ and $b=\LL v,J\gamma'(\ell)\RL$, then $v$ has a decomposition
    \beq
        v=a \gamma'(\ell)+bJ\gamma'(\ell)+v_0,
    \eeq
where $\LL v_0,\gamma'(\ell)\RL=0$ and $\LL v_0,J\gamma'(\ell)\RL=0$. Moreover,
    \beq
        Jv=-b\gamma'(\ell)+aJ\gamma'(\ell)+Jv_0.
    \eeq
    Furthermore, if we choose two normal vectors
    \beq
        v_1=bJ\gamma'(\ell)+v_0\qtq{and} v_2=aJ\gamma'(\ell)+Jv_0,
    \eeq
    then one has
    \beq
        \left(\mathrm{Hess}\,r\right)(v,v)+\left(\mathrm{Hess}\,r\right)(Jv,J v)
        = \left(\mathrm{Hess}\,r\right)(v_1,v_1)+\left(\mathrm{Hess}\,r\right)(v_2,v_2).
    \eeq
    Consider two variational vector fields
    \beq
        U_1=\tilde V-\LL \tilde V,\gamma'\RL\gamma' \qtq{and}
        U_2=J\tilde V-\LL J\tilde V,\gamma'\RL\gamma'.
    \eeq
    Since $V(0)=0$, one has $U_1(0)=U_2(0)=0$. Moreover, $\LL U_1,\gamma'\RL=\LL U_2,\gamma'\RL\equiv 0$,
    \beq
        U_1(\ell)= v-\LL v,\gamma'(\ell)\RL\gamma'(\ell) = v- a\gamma'(\ell)=v_1,
    \eeq
    \beq
        U_2(\ell)= Jv-\LL Jv,\gamma'(\ell)\RL\gamma'(\ell) = Jv+b\gamma'(\ell)=v_2.
    \eeq
    By using the index form theorem in Riemannian geometry, one deduces that
    \begin{eqnarray}
        2\left(\p\bp r\right)\left(X,\bar X\right)
        &=& \left(\mathrm{Hess}\,r\right)(v,v)+\left(\mathrm{Hess}\,r\right)(Jv,J v)  \nonumber\\
        &=& \left(\mathrm{Hess}\,r\right)(v_1,v_1)+\left(\mathrm{Hess}\,r\right)(v_2,v_2) \nonumber\\
        &\leq&  I_\gamma(U_1,U_1)+   I_\gamma(U_2,U_2).\label{index form}
    \end{eqnarray}
    Since $U_1'=\tilde V'-\LL \tilde V',\gamma'\RL\gamma' $, one has
    \begin{eqnarray}
        I_\gamma(U_1,U_1)
        &=& \int_0^\ell |U_1'|^2-R(U_1,\gamma',\gamma',U_1)\,dt  \nonumber\\
        &=& \int_0^\ell |\tilde V'|^2-\LL\tilde V',\gamma'\RL^2-R(\tilde V,\gamma',\gamma',\tilde V)\,dt \nonumber\\
        &=& I_\gamma(\tilde V,\tilde V)-\int_0^\ell \LL\tilde V',\gamma'\RL^2\,dt.
    \end{eqnarray}
    Similarly, one can derive
    \beq
        I_\gamma(U_2,U_2)= I_\gamma(J\tilde V,J\tilde V)-\int_0^\ell \LL  J\tilde V ',\gamma'\RL^2\,dt.
    \eeq
    Therefore, one deduces that
    \begin{eqnarray}
        2\left(\p\bp r\right)\left(X,\bar X\right)
        &\leq&
            I_\gamma(\tilde V,\tilde V)+I_\gamma(J\tilde V,J\tilde V)
        - \int_0^\ell \LL \tilde V',\gamma'\RL^2 dt- \int_0^\ell \LL J\tilde V',\gamma'\RL^2 dt\nonumber\\
        &=&  I_\gamma(\tilde V,\tilde V)+  I_\gamma(J\tilde V,J\tilde V) -\int_0^\ell \left|\LL V',E_\gamma\RL\right|^2dt,
        \label{d d bar r key estimate}
    \end{eqnarray}
    where   we use the elementary fact that
    \beq
    \LL \tilde V',\gamma'\RL^2\ +\LL J\tilde V',\gamma'\RL^2=\left|\LL V',E_\gamma\RL\right|^2.
    \eeq
    On the other hand,  since $$R\left(E_\gamma,\bar{E_\gamma},V,\bar V\right)=R\left(\gamma',\tilde V,\tilde V,\gamma'\right) + R\left(\gamma',J\tilde V,J\tilde V,\gamma'\right),$$ one obtains
    \begin{eqnarray}
        \cX_\gamma\left(V,\bar V\right)
        &=& \int_0^\ell |V'|^2 -\frac{1}{2}R\left(E_\gamma,\bar{E_\gamma},V,\bar V\right)dt \nonumber\\
        &=& \int_0^\ell \frac{1}{2}|\tilde V'|^2 +\frac{1}{2}|J\tilde V'|^2-\frac{1}{2}R\left(\gamma',\tilde V,\tilde V,\gamma'\right) -\frac{1}{2} R\left(\gamma',J\tilde V,J\tilde V,\gamma'\right)dt \nonumber\\
        &=& \frac{1}{2} I_\gamma(\tilde V,\tilde V)+ \frac{1}{2} I_\gamma(J\tilde V,J\tilde V) ,\label{index form calculation}
    \end{eqnarray}
    By \eqref{d d bar r key estimate} and \eqref{index form calculation}, one establishes the estimate in \eqref{new index form estimate}.
    Assuming the identity in \eqref{new index form estimate} holds, it follows that the identity in \eqref{index form} is satisfied. Consequently, both $U_1$ and
    $U_2$   are Jacobi fields. This implies that $V^\perp=V-\LL V,\gamma'\RL\gamma'$ is a complex  Jacobi field. The converse statement is obvious.
\end{proof}

\noindent  The following result is standard and will be invoked repeatedly in subsequent arguments.
\btheorem \label{thm from Jacobi to curvature}
Let $(M,\om_g)$ be a complete K\"ahler manifold, $p\in M$ and $U=M\backslash\mathrm{cut}{(p)}$.
    Then the following are equivalent.
    \bd
    \item $(M,\om_g)$ has constant holomorphic bisectional $c\in\R$.
    \item  Let $\gamma:[0,\ell]\to U$ be a unit speed geodesic with $\gamma(0)=p$.
        Every Jacobi field along $\gamma$ with $J(0)=0$ and $\LL J,\gamma'\RL\equiv 0$ is of the form
        \beq
            J(t)=a \operatorname{sn}_{c/2}(t)E(t)+ b\operatorname{sn}_{2c}(t)J\gamma'(t)\label{canonical Jacobi field}
        \eeq
        where $E(t)$ is any parallel vector field along $\gamma$ with $\LL E(t),\gamma'(t)\RL=\LL E(t),J\gamma'(t)\RL\equiv 0$ and $|E(t)|\equiv 1$.
    \ed
\etheorem

\vskip \baselineskip

\section{Comparison theorems for   symmetrized curvature operators}

 In this section we prove Theorem \ref{thm symmetric curvature operator comparison},  Theorem \ref{symmetric laplacian positive} and Theorem \ref{thm diam}. Let $A$ be a $k$-semi-positive Hermitian $ n \times n $ matrix and $ \lambda_1 \leq \cdots \leq \lambda_n $ be eigenvalues of $ A $.   Suppose that
$ \{e_i\}_{i=1}^n $ is an orthonormal basis of $ \C^{n} $. One deduces that
\beq  \sum_{s=1}^k \langle Ae_{i_s}, e_{i_s} \rangle \geq \lambda_1 + \cdots + \lambda_k\geq 0, \label{partialtrace} \eeq
for any $1\leq i_1<\cdots<i_k\leq n$.

\blemma\label{keylinear}  Fix $c\in \R$ and $k<n$.  If
$\sR-2c\cdot\id$ is $k$-semi-positive at point $p\in M$, then for
any orthonormal basis $E_1,\cdots,E_{n}$ of $ T^{1,0}_pM$ and any
$\alpha\geq \frac{2(k-1)}{n-1}$, one has \beq
R\left(E_n,\bar{E_n},E_n,\bar{E_n}\right) +\alpha
\sum_{i=1}^{n-1}R\left(E_n,\bar{E_n},E_i,\bar{E_i}\right)\geq
2c+\alpha(n-1)c. \label{mixed estimate}\eeq \elemma

\begin{proof}
    Consider the following orthonormal vectors in $\mathrm{Sym}^2 T^{1,0}_pM$:
    \beq
    V_n=E_n\otimes E_n,\quad V_i=\frac{E_n\otimes E_i+E_i\otimes E_n}{\sqrt{2}},
    \eeq
    where $1\leq i\leq n-1$. Since $\sR-2c\cdot\id$ is $k$-semi-positive, for any subset $I\subset \{1,\cdots,n-1\}$ with $|I|=k-1$, we have
    \beq
    \sR\left(V_n,\bar{V_n}\right)+\sum_{i\in I} \sR\left(V_i,\bar{V_i}\right)\geq 2ck.
    \eeq
    Summing over all such subsets $I$ and taking average, we conclude that
    \beq
    \sR\left(V_n,\bar{V_n}\right)+\frac{k-1}{n-1}\sum_{i=1}^{n-1} \sR\left(V_i,\bar{V_i}\right)\geq 2ck.
    \eeq
    On the other hand, a direct calculation shows that
    \beq
    \sR\left(V_n,\bar{V_n}\right) = R\left(E_n,\bar{E_n},E_n,\bar{E_n}\right),
    \eeq
    and for $1\leq i\leq n-1$,
    \beq
    \sR\left(V_i,\bar{V_i}\right)=2\sR\left(E_n\otimes E_i,\bar{E_n}\otimes \bar{E_i}\right)=2R\left(E_n,\bar{E_n},E_i,\bar{E_i}\right).
    \eeq
    Thus we obtain
    \beq
    R\left(E_n,\bar{E_n},E_n,\bar{E_n}\right)+\frac{2(k-1)}{n-1}\sum_{i=1}^{n-1}R\left(E_n,\bar{E_n},E_i,\bar{E_i}\right)\geq 2ck.
    \label{average 1}
    \eeq
    Moreover, for any subset $J\subset \{1,\cdots,n-1\}$ with $|J|=k$, we have
    \beq
    \sum_{i\in J} \sR\left(V_i,\bar{V_i}\right)\geq 2ck,
    \eeq
    and we deduce that
    \beq
    \sum_{i=1}^{n-1}R\left(E_n,\bar{E_n},E_i,\bar{E_i}\right)\geq (n-1)c.
    \label{average 2}
    \eeq
    By \eqref{average 1} and \eqref{average 2}, we obtain
    \be && R\left(E_n,\bar{E_n},E_n,\bar{E_n}\right)
    +\alpha \sum_{i=1}^{n-1}R\left(E_n,\bar{E_n},E_i,\bar{E_i}\right)\\
    &=&  R\left(E_n,\bar{E_n},E_n,\bar{E_n}\right)
    +\frac{2(k-1)}{n-1} \sum_{i=1}^{n-1}R\left(E_n,\bar{E_n},E_i,\bar{E_i}\right)\\&&+\left(\alpha-\frac{2(k-1)}{n-1}\right) \sum_{i=1}^{n-1}R\left(E_n,\bar{E_n},E_i,\bar{E_i}\right)\\
    &\geq & 2ck+\left(\alpha-\frac{2(k-1)}{n-1}\right)(n-1)c=2c+\alpha(n-1)c.
    \ee
    This completes the proof.
\end{proof}

\bcorollary\label{symmetry to ric} If $\sR-2c\cdot \mathrm{id}$ is
$k$-semi-positive for some $k\leq (n+1)/2$, then \beq
\mathrm{Ric}(\omega_g)\geq (n+1)c\om_g. \eeq \ecorollary

\noindent The following is Theorem \ref{thm symmetric curvature
operator comparison}: \btheorem\label{symmetric laplacian negative}
Let $(M,\om_g)$ be a complete K\"ahler manifold of complex dimension
$n$ and $r$ be  the distance function from a given point $p\in M$.
If $\sR-2c\cdot \mathrm{id}$ is $k$-semi-positive for some $c<0$ and $k\leq
(n+1)/2$, then \beq \Delta
r\leq2(n-1)\frac{\mathrm{sn}_{c/2}'(r)}{\mathrm{sn}_{c/2}(r)}+\frac{\mathrm{sn}_{2c}'(r)}{\mathrm{sn}_{2c}(r)}\label{Laplacian
comparison 2} \eeq on $M\setminus \mathrm{cut}(p)\cup \{p\}$.
Moreover, if the identity in \eqref{Laplacian comparison 2} holds on
$M\setminus\mathrm{cut}(p)\cup \{p\}$, then the universal cover of
$(M,\omega_g)$ is isometrically biholomorphic to the hyperbolic
space $M_c$. \etheorem \bproof     For a point $x\in
M\setminus\mathrm{cut}(p)\cup \{p\}$, let $\gamma\colon [0,\ell]\to
M$ be the unit-speed minimal geodesic joining $p$ and $x$, and
$E_1(t),\cdots,E_n(t)\in\Gamma([0,\ell],\gamma^*T^{1,0}M)$ be
orthonormal parallel fields along $\gamma$ such that
$E_n=\frac{1}{\sqrt{2}}\left(\gamma'-\sq J\gamma'\right)$. We define
vector fields \beq
V_n(t)=\frac{\mathrm{sn}_{2c}(t)}{\mathrm{sn}_{2c}(\ell)}
E_n(t)\qtq{and}
V_i(t)=\frac{\mathrm{sn}_{c/2}(t)}{\mathrm{sn}_{c/2}(\ell)} E_i(t)
\eeq where $1\leq i\leq n-1$. By Theorem \ref{thm new index}, one
has \be \left(\p\bp r\right)\left(E_n(\ell),\bar{E_n(\ell)}\right)
&\leq&  \cX_{\gamma}\left(V_n,\bar {V_n}\right)-\frac{1}{2}\int_0^\ell \left|\LL V_n',E_n\RL\right|^2 dt\nonumber\\
&=&  \frac{1}{2} \int_0^\ell  \frac{\mathrm{sn}'_{2c}(t)^2}{\mathrm{sn}_{2c}(\ell)^2}\,dt
-\frac{1}{2}\int_0^\ell  \frac{\mathrm{sn}_{2c}(t)^2}{\mathrm{sn}_{2c}(\ell)^2} R\left(E_n,\bar{E_n},E_n,\bar{E_n}\right) dt.
\ee
Moreover, for $1\leq i\leq n-1$,
\be
\left(\p\bp r\right)\left(E_i(\ell),\bar{E_i(\ell)}\right)
&\leq& \cX_{\gamma}\left(V_i,\bar {V_i}\right)\nonumber\\
&=&\int_0^\ell  \frac{\mathrm{sn}'_{c/2}(t)^2}{\mathrm{sn}_{c/2}(\ell)^2}\,dt
-\frac{1}{2}\int_0^\ell \frac{\mathrm{sn}_{c/2}(t)^2}{\mathrm{sn}_{c/2}(\ell)^2}R\left(E_n,\bar{E_n},E_i,\bar{E_i}\right) dt.
\ee
It follows that
\begin{eqnarray}
\Delta r(x)
&=& 2\sum_{i=1}^n\left(\p\bp r\right)\left(E_i(\ell),\bar{E_i(\ell)}\right) \nonumber\\
&\leq&   \int_0^\ell \left\{ \frac{\mathrm{sn}'_{2c}(t)^2}{\mathrm{sn}_{2c}(\ell)^2}+2(n-1)\frac{\mathrm{sn}'_{c/2}(t)^2}{\mathrm{sn}_{c/2}(\ell)^2}\right\}dt \label{Laplacian estimate 2} \\
&&\nonumber-\int_0^\ell  \left\{\frac{\mathrm{sn}_{2c}(t)^2}{\mathrm{sn}_{2c}(\ell)^2} R\left(E_n,\bar{E_n},E_n,\bar{E_n}\right)
+\frac{\mathrm{sn}_{c/2}(t)^2}{\mathrm{sn}_{c/2}(\ell)^2}\sum_{i=1}^{n-1}R\left(E_n,\bar{E_n},E_i,\bar{E_i}\right) \right\}dt.
\end{eqnarray}
For $t\in[0,\ell]$, we set
\beq \alpha(t):=\frac{\mathrm{sn}_{c/2}^2(t)}{\mathrm{sn}_{c/2}^2(\ell)}\cdot \left(\frac{\mathrm{sn}_{2c}^2(t)}{\mathrm{sn}_{2c}^2(\ell)}\right)^{-1}.\eeq
A simple calculation shows
$$\alpha(t)=\left(\frac{e^{\sqrt{-c/2}\cdot t}+e^{-\sqrt{-c/2}\cdot t}}{e^{\sqrt{-c/2}\cdot \ell}+e^{-\sqrt{-c/2}\cdot \ell}}\right)^{-2}\geq \alpha(\ell)= 1. $$
Since $k\leq (n+1)/2$, \beq \alpha(t)\geq 1 \geq \frac{2(k-1)}{n-1}. \eeq
We  apply Lemma \ref{keylinear} to this $\alpha(t)$
and obtain
\be
&&\frac{\mathrm{sn}_{2c}(t)^2}{\mathrm{sn}_{2c}(\ell)^2} R\left(E_n,\bar{E_n},E_n,\bar{E_n}\right)
+\frac{\mathrm{sn}_{c/2}(t)^2}{\mathrm{sn}_{c/2}(\ell)^2}\sum_{i=1}^{n-1}R\left(E_n,\bar{E_n},E_i,\bar{E_i}\right)\\
&=&\frac{\mathrm{sn}_{2c}(t)^2}{\mathrm{sn}_{2c}(\ell)^2}\left(R\left(E_n,\bar{E_n},E_n,\bar{E_n}\right)
+\alpha(t) \sum_{i=1}^{n-1}R\left(E_n,\bar{E_n},E_i,\bar{E_i}\right)\right)\\
&\geq& 2c \frac{\mathrm{sn}_{2c}(t)^2}{\mathrm{sn}_{2c}(\ell)^2} + (n-1)c \frac{\mathrm{sn}_{c/2}(t)^2}{\mathrm{sn}_{c/2}(\ell)^2}.\label{mixed term 2}
\ee
One deduces that
\be
\Delta r(x)
&\leq&   \int_0^\ell \left\{ \frac{\mathrm{sn}'_{2c}(t)^2}{\mathrm{sn}_{2c}(\ell)^2}+2(n-1)\frac{\mathrm{sn}'_{c/2}(t)^2}{\mathrm{sn}_{c/2}(\ell)^2}\right\}dt \nonumber\\
&&-\int_0^\ell  \left\{2c \frac{\mathrm{sn}_{2c}(t)^2}{\mathrm{sn}_{2c}(\ell)^2} + (n-1)c \frac{\mathrm{sn}_{c/2}(t)^2}{\mathrm{sn}_{c/2}(\ell)^2}\right\}dt\nonumber\\
&=& 2(n-1)\frac{\mathrm{sn}_{c/2}'(\ell)}{\mathrm{sn}_{c/2}(\ell)}+\frac{\mathrm{sn}_{2c}'(\ell)}{\mathrm{sn}_{2c}(\ell)}.
\ee
This completes the proof of \eqref{Laplacian comparison 2}.

Moreover, if the identity in \eqref{Laplacian comparison 2} holds on $M\setminus\mathrm{cut}(p)\cup \{p\}$,
then the identity in \eqref{Laplacian estimate 2} holds.
Then by Theorem \ref{thm new index} and  Theorem \ref{thm from Jacobi to curvature}, $(M,\om_g)$ has constant holomorphic bisectional curvature $c<0$,
and so the universal cover of $(M,\omega_g)$ is isometrically biholomorphic to the hyperbolic space $M_c$.
\eproof

\noindent We prove Theorem \ref{symmetric laplacian positive}:

\btheorem\label{main2} Let $(M,\om_g)$ be a complete K\"ahler
manifold of complex dimension $n$ and $r$ be  the distance function
from a given point $p\in M$. If $\sR-2c\cdot \mathrm{id}$ is
$k$-semi-positive for some $c>0$ and $k< (n+1)/2$, then \beq \Delta
r\leq2(n-1)\frac{\mathrm{sn}_{c/2}'(r)}{\mathrm{sn}_{c/2}(r)}+\frac{\mathrm{sn}_{2c}'(r)}{\mathrm{sn}_{2c}(r)}\label{Laplacian
comparison 3} \eeq on some metric ball
$B(p,C(k,n)/\sqrt{c})\setminus \mathrm{cut}(p)\cup \{p\}$. Moreover
if the identity \eqref{Laplacian comparison 3} holds for all such
$x$, then $\mathrm{HBSC}\equiv c$ on $B(p,C(k,n)/\sqrt{c})$.
\etheorem

\bproof Using the same setup  as in the proof of Theorem \ref{symmetric laplacian negative}, we have
\begin{eqnarray}
\Delta r(x)&\leq &
 \int_0^\ell \left\{ \frac{\mathrm{sn}'_{2c}(t)^2}{\mathrm{sn}_{2c}(\ell)^2}+2(n-1)\frac{\mathrm{sn}'_{c/2}(t)^2}{\mathrm{sn}_{c/2}(\ell)^2}\right\}dt \label{Laplacian estimate 3} \\
&&\nonumber-\int_0^\ell  \left\{\frac{\mathrm{sn}_{2c}(t)^2}{\mathrm{sn}_{2c}(\ell)^2} R\left(E_n,\bar{E_n},E_n,\bar{E_n}\right)
+\frac{\mathrm{sn}_{c/2}(t)^2}{\mathrm{sn}_{c/2}(\ell)^2}\sum_{i=1}^{n-1}R\left(E_n,\bar{E_n},E_i,\bar{E_i}\right) \right\}dt.
\end{eqnarray}
For $t\in[0,\ell]$, we set
\beq \alpha(t):=\frac{\mathrm{sn}_{c/2}^2(t)}{\mathrm{sn}_{c/2}^2(\ell)}\cdot \left(\frac{\mathrm{sn}_{2c}^2(t)}{\mathrm{sn}_{2c}^2(\ell)}\right)^{-1}=\left(\frac{\cos\left(\sqrt{\frac{c}{2}}\ell\right)}{\cos\left(\sqrt{\frac{c}{2}}t\right)}\right)^2\geq \cos^{2}\left(\sqrt{\frac{c}{2}}\ell\right).
\eeq
Hence if $k<(n+1)/2$ and $\alpha(t)\geq \frac{2(k-1)}{n-1}$, we can apply Lemma \ref{keylinear} and obtain
\be
&&\frac{\mathrm{sn}_{2c}(t)^2}{\mathrm{sn}_{2c}(\ell)^2} R\left(E_n,\bar{E_n},E_n,\bar{E_n}\right)
+\frac{\mathrm{sn}_{c/2}(t)^2}{\mathrm{sn}_{c/2}(\ell)^2}\sum_{i=1}^{n-1}R\left(E_n,\bar{E_n},E_i,\bar{E_i}\right)\\
&\geq& 2c \frac{\mathrm{sn}_{2c}(t)^2}{\mathrm{sn}_{2c}(\ell)^2} + (n-1)c \frac{\mathrm{sn}_{c/2}(t)^2}{\mathrm{sn}_{c/2}(\ell)^2}.\label{mixed term 3}
\ee
Note that
$ \alpha(t)\geq \frac{2(k-1)}{n-1}$ for all $t\in[0,\ell]$ if and only if
$\cos^{2}\left(\sqrt{\frac{c}{2}}\ell\right)\geq \frac{2(k-1)}{n-1}$. Hence,  $ \alpha(t)\geq \frac{2(k-1)}{n-1}$ holds for all $t\in[0,\ell]$ if $\ell\leq  C(k,n)/\sqrt{c}$
for some constant $C(k,n)$.  By using similar arguments as in the proof of Theorem \ref{symmetric laplacian negative} we obtain
$$
\Delta r\leq2(n-1)\frac{\mathrm{sn}_{c/2}'(r)}{\mathrm{sn}_{c/2}(r)}+\frac{\mathrm{sn}_{2c}'(r)}{\mathrm{sn}_{2c}(r)}
$$
on $B(p,C(k,n))\setminus \mathrm{cut}(p)\cup \{p\}$. The rigidity property is derived from the proofs of Theorem \ref{thm new index} and Theorem \ref{thm from Jacobi to curvature}.
\eproof

\noindent The following is Theorem \ref{thm diam}: \btheorem Let
$(M, \omega_g)$ be a complete K\"ahler manifold of complex dimension
$n$. If $\sR-2c\cdot\id$ is $k$-semi-positive for some $c>0$ and
$k\leq n$, then $M$ is compact and \beq \mathrm{diam}(M,\om_g)\leq
\frac{\pi}{\sqrt{\nu}}\qtq{where}\nu=\frac{2kc}{4k-3}. \eeq
Moreover, the fundamental group $\pi_1(M)$ of $M$ is trivial.

\etheorem

\begin{proof}
    Let $d_0$ denote the diameter of $(M,\omega_g)$.
    If $d_0 > \frac{\pi}{\sqrt{\nu}}$, then there exist points $p, q \in M$ such that $d=d(p,q) > \frac{\pi}{\sqrt{\nu}}$.
    Let $\gamma \colon [0,d] \to M$ be a unit-speed minimal geodesic joining $p$ and $q$. Let $E_1(t), \cdots, E_n(t) \in \Gamma([0,d], \gamma^*T^{1,0}M)$ be orthonormal parallel vector fields along $\gamma$ such that $E_n = \frac{1}{\sqrt{2}}(\gamma' - \sqrt{-1} J\gamma').$
    Fix $\ell \in \left( 0, \frac{\pi}{\sqrt{\nu}}\right)\subset (0,d)$, and define  vector fields
    \beq
    V_i(t) = \frac{\mathrm{sn}_{\nu}(t)}{\mathrm{sn}_{\nu}(\ell)} E_i(t), \quad 0 \leq t \leq \ell,
    \eeq
    along $\gamma|_{[0,\ell]}$ for $1 \leq i \leq n-1$.
    Since $\gamma(\ell)\not\in\mathrm{cut}(p)$, by Theorem \ref{thm new index}, one has
    \be
    \left(\p\bp r\right)\left(E_n(\ell),\bar{E_n(\ell)}\right)
    &\leq&  \cX_{\gamma|_{[0,\ell]}}\left(V_n,\bar {V_n}\right)-\frac{1}{2}\int_0^\ell \left|\LL V_n',E_n\RL\right|^2 dt\nonumber\\
    &=&  \frac{1}{2} \int_0^\ell  \frac{\mathrm{sn}'_{\nu}(t)^2}{\mathrm{sn}_{\nu}(\ell)^2}\,dt
    -\frac{1}{2}\int_0^\ell  \frac{\mathrm{sn}_{\nu}(t)^2}{\mathrm{sn}_{\nu}(\ell)^2} R\left(E_n,\bar{E_n},E_n,\bar{E_n}\right) dt.
    \ee
    Moreover, for $1\leq i\leq n-1$,
    \be
    \left(\p\bp r\right)\left(E_i(\ell),\bar{E_i(\ell)}\right)
    &\leq& \cX_{\gamma|_{[0,\ell]}}\left(V_i,\bar {V_i}\right)\nonumber\\
    &=&\int_0^\ell  \frac{\mathrm{sn}'_{\nu}(t)^2}{\mathrm{sn}_{\nu}(\ell)^2}\,dt
    -\frac{1}{2}\int_0^\ell \frac{\mathrm{sn}_{\nu}(t)^2}{\mathrm{sn}_{\nu}(\ell)^2}R\left(E_n,\bar{E_n},E_i,\bar{E_i}\right) dt.
    \ee
Consider the following combination:
    \begin{eqnarray}
    &&\frac{1}{2} \left(\p\bp r\right)\left(E_n(\ell),\bar{E_n(\ell)}\right) +
    \sum_{i=1}^{k-1} \left(\p\bp r\right)\left(E_i(\ell),\bar{E_i(\ell)}\right)\nonumber\\
    &\leq & \frac{4k-3}{4} \int_0^\ell  \frac{\mathrm{sn}'_{\nu}(t)^2}{\mathrm{sn}_{\nu}(\ell)^2}\,dt \\
    &&-\frac{1}{4}\int_0^\ell  \frac{\mathrm{sn}_{\nu}(t)^2}{\mathrm{sn}_{\nu}(\ell)^2} \left\{R\left(E_n,\bar{E_n},E_n,\bar{E_n}\right)+2\sum_{i=1}^{k-1}R\left(E_n,\bar{E_n},E_i,\bar{E_i}\right)\right\}dt.\nonumber
    \end{eqnarray}
    By using the argument in the proof of Lemma~\ref{keylinear}, we obtain
    \beq \sR\left(V_n,\bar{V_n}\right)+\sum_{i=1}^{k-1} \sR\left(V_i,\bar{V_i}\right)=
    R\left(E_n,\bar{E}_n,E_n,\bar{E}_n\right) + 2\sum_{i=1}^{k-1} R\left(E_n,\bar{E}_n,E_i,\bar{E}_i\right) \geq 2kc.
    \eeq
    Therefore,
    \begin{eqnarray}
    &&\frac{1}{2} \left(\p\bp r\right)\left(E_n(\ell),\bar{E_n(\ell)}\right) +
    \sum_{i=1}^{k-1} \left(\p\bp r\right)\left(E_i(\ell),\bar{E_i(\ell)}\right) \label{diam 1} \\
    &\leq& \frac{4k-3}{4} \int_0^\ell  \frac{\mathrm{sn}'_{\nu}(t)^2}{\mathrm{sn}_{\nu}(\ell)^2}\,dt -\frac{kc}{2}\int_0^\ell  \frac{\mathrm{sn}_{\nu}(t)^2}{\mathrm{sn}_{\nu}(\ell)^2}  dt
    = \frac{4k-3}{4}\frac{\mathrm{sn}'_{\nu}(\ell)}{\mathrm{sn}_{\nu}(\ell)}.\nonumber
    \end{eqnarray}
    Since $d>\frac{\pi}{\sqrt{\nu}}$, the distance function $r$ is smooth at $\gamma\left(\frac{\pi}{\sqrt{\nu}}\right)$, and so there exists $\delta>0$ such that the left-hand side of \eqref{diam 1} is uniformly bounded from below for all $\ell\in \left(\frac{\pi}{\sqrt{\nu}}-\delta,\frac{\pi}{\sqrt{\nu}}\right)$. However, the right-hand side has the property that
    \beq
    \lim_{t\nearrow \frac{\pi}{\sqrt{\nu}}}\frac{\mathrm{sn}_{\nu}'(t)}{\mathrm{sn}_{\nu}(\ell)}=-\infty.
    \eeq
This is a contradiction.

    Let $\pi:(\tilde M,\tilde \om)\to (M,\om_g)$ be the universal cover. Then $(\tilde M,\tilde \om)$ has the same curvature property and so it is also compact.
Therefore,  $\pi$ must be a finite cover. This implies that $\pi_1(M)$ is finite.
Furthermore, the finiteness of $\pi_1(M)$ indicates $H^1(M, \mathbb{C}) = 0$, and consequently, $\dim H^{0,1}_{\bp}(M,\C)=0.$
On the other hand, by \cite[Theorem~1.5]{WY25+}, we obtain
\beq
H^{0,i}_{\bp}(M,\C)=0
\eeq
for $2\leq i\leq n$. Therefore, the holomorphic Euler characteristic
\beq
\chi(M,\cO_M)
=\sum_{i=0}^n(-1)^i \dim H^{0,i}_{\bp}(M,\C)
= \dim H^{0,0}_{\bp}(M,\C) =1.
\eeq
Since the universal cover $(\tilde M,\tilde \om)$ has the same curvature property, we also have $
\chi\left(\tilde M,\cO_{\tilde M}\right) =1$.
The Riemann-Roch theorem asserts that
\beq
\chi\left(\tilde M,\cO_{\tilde M}\right)
=\left|\pi_1(M)\right|\cdot \chi(M,\cO_M).
\eeq
Therefore, $\left|\pi_1(M)\right|=1$, i.e., the fundamental group $\pi_1(M)$ is trivial.\end{proof}

\noindent The following result is of particular interest. \btheorem
Let $(M,\om_g)$ be a complete K\"ahler manifold of complex dimension
$n$ and $r$ be the distance function from a given point $p\in M$. If
$\sR-2c\cdot \mathrm{id}$ is $k$-semi-positive for some $c\in\R$ and $k<n$, then at
a point $x\in M\setminus \mathrm{cut}(p)\cup\{p\}$, one has \beq
    \sum_{i=1}^k\left(\p\bp r\right)\left(X_i,\bar{X_i}\right)\leq
    k\frac{\mathrm{sn}_{c/2}'(r)}{\mathrm{sn}_{c/2}(r)}
    \label{khessiancomparison}
\eeq
where $E_r|_x,X_1,\cdots,X_{k}$ are orthonormal in $T_x^{1,0}M$ and
$$E_r=\frac{\nabla r-\sq J\nabla r}{\sqrt{2}}.$$
\etheorem

\bproof
    We use the same notation as in the proof of Theorem \ref{symmetric laplacian negative}, and assume $X_i=E_i(\ell)$ for $1\leq i\leq k$, $E_n=E_r$. For every $1\leq i\leq k$, we  have
    \beq
        \left(\p\bp r\right)\left(X_i,\bar{X_i}\right)
        \leq \int_0^\ell  \frac{\mathrm{sn}'_{c/2}(t)^2}{\mathrm{sn}_{c/2}(\ell)^2}\,dt
        -\frac{1}{2}\int_0^\ell \frac{\mathrm{sn}_{c/2}(t)^2}{\mathrm{sn}_{c/2}(\ell)^2}R\left(E_n,\bar{E_n},E_i,\bar{E_i}\right) dt.
    \eeq
    Following the argument in the proof of Lemma~\ref{keylinear}, one obtains
    \beq
    \sum_{i=1}^{k-1} \sR\left(V_i,\bar{V_i}\right)=2    \sum_{i=1}^k R\left(E_n,\bar{E_n},E_i,\bar{E_i}\right)\geq 2ck.
    \eeq
    Therefore, one has
    \begin{eqnarray}
        \sum_{i=1}^k\left(\p\bp r\right)\left(X_i,\bar{X_i}\right)
        &\leq& k\int_0^\ell  \frac{\mathrm{sn}'_{c/2}(t)^2}{\mathrm{sn}_{c/2}(\ell)^2}\,dt
        -\frac{1}{2}\int_0^\ell \frac{\mathrm{sn}_{c/2}(t)^2}{\mathrm{sn}_{c/2}(\ell)^2} \sum_{i=1}^k R\left(E_n,\bar{E_n},E_i,\bar{E_i}\right) dt\nonumber\\
        &\leq& k\int_0^\ell  \frac{\mathrm{sn}'_{c/2}(t)^2}{\mathrm{sn}_{c/2}(\ell)^2}\,dt
        -\frac{ck}{2}\int_0^\ell \frac{\mathrm{sn}_{c/2}(t)^2}{\mathrm{sn}_{c/2}(\ell)^2} dt
        = k\frac{\mathrm{sn}_{c/2}'(\ell)}{\mathrm{sn}_{c/2}(\ell)}.
    \end{eqnarray}
    This completes the proof.
\eproof

\vskip 2\baselineskip

\section{Comparison theorems for Ricci curvature and HSC}

In this section we prove Theorem \ref{thm Laplacian Ric and HSC} and Theorem \ref{thm volume Ric and HSC}.

\blemma\label{keylinear2} Let $(M,\omega_g)$ be a K\"ahler manifold of complex dimension $n$.  If $$\mathrm{HSC}\geq 2c \qtq{and} \mathrm{Ric}\geq c(n+1)\om_g$$ at  point $p\in M$ for some $c\in \R$, then for any orthonormal basis $E_1,\cdots,E_{n}\in T^{1,0}_pM$ and any real number $\alpha\in [0,1]$, one has
\beq
R\left(E_n,\bar{E_n},E_n,\bar{E_n}\right)
+\alpha \sum_{i=1}^{n-1}R\left(E_n,\bar{E_n},E_i,\bar{E_i}\right)\geq 2c+\alpha(n-1)c. \label{mixed estimate2}\eeq
\elemma

\bproof  It is easy to see that
\be &&R\left(E_n,\bar{E_n},E_n,\bar{E_n}\right)
+\alpha \sum_{i=1}^{n-1}R\left(E_n,\bar{E_n},E_i,\bar{E_i}\right)\\&=&(1-\alpha)R\left(E_n,\bar{E_n},E_n,\bar{E_n}\right)+\alpha \mathrm{Ric}(E_n,E_n)\\&\geq&(1-\alpha)\cdot 2c+\alpha(n+1)c=2c+\alpha(n-1)c.  \ee
This completes the proof.
\eproof

\noindent We restate Theorem \ref{thm Laplacian Ric and HSC} here
for readers' convenience. \btheorem Let $(M, \omega_g)$ be a
complete K\"ahler manifold of complex dimension $n$. If
$\mathrm{HSC}\geq 2c$ and $\mathrm{Ric}\geq c(n+1)\om_g$ for some
$c>0$, then the distance function $r$ from a  given point $p\in M$
satisfies \beq \Delta r\leq
2(n-1)\frac{\mathrm{sn}_{c/2}'(r)}{\mathrm{sn}_{c/2}(r)}+\frac{\mathrm{sn}_{2c}'(r)}{\mathrm{sn}_{2c}(r)}\label{Laplacian
comparison 1} \eeq on $M\setminus\mathrm{cut}(p)\cup \{p\}$.
Moreover, if the identity in \eqref{Laplacian comparison 1} holds on
$M\setminus\mathrm{cut}(p)\cup \{p\}$, then $(M,\omega_g)$ is
isometrically biholomorphic to the projective space $M_c$. \etheorem

\bproof Using  the same notation as in the proof of Theorem \ref{symmetric laplacian negative}, we have
\begin{eqnarray}
\Delta r(x)&\leq &
\int_0^\ell \left\{ \frac{\mathrm{sn}'_{2c}(t)^2}{\mathrm{sn}_{2c}(\ell)^2}+2(n-1)\frac{\mathrm{sn}'_{c/2}(t)^2}{\mathrm{sn}_{c/2}(\ell)^2}\right\}dt \label{Laplacian estimate 3} \\
&&\nonumber-\int_0^\ell  \left\{\frac{\mathrm{sn}_{2c}(t)^2}{\mathrm{sn}_{2c}(\ell)^2} R\left(E_n,\bar{E_n},E_n,\bar{E_n}\right)
+\frac{\mathrm{sn}_{c/2}(t)^2}{\mathrm{sn}_{c/2}(\ell)^2}\sum_{i=1}^{n-1}R\left(E_n,\bar{E_n},E_i,\bar{E_i}\right) \right\}dt.
\end{eqnarray}
For $t\in[0,\ell]$, we set
\beq \alpha(t):=\frac{\mathrm{sn}_{c/2}^2(t)}{\mathrm{sn}_{c/2}^2(\ell)}\cdot \left(\frac{\mathrm{sn}_{2c}^2(t)}{\mathrm{sn}_{2c}^2(\ell)}\right)^{-1}=\left(\frac{\cos\left(\sqrt{\frac{c}{2}}\ell\right)}{\cos\left(\sqrt{\frac{c}{2}}t\right)}\right)^2.
\eeq
Since $\mathrm{HSC}\geq 2c$, by \cite{Tsu57}, $\mathrm{diam}(M,\omega_g)\leq \frac{\pi}{\sqrt{2c}}$ and so $\ell< \frac{\pi}{\sqrt{2c}} $. We conclude that $\cos\left(\sqrt{\frac{c}{2}}t\right)$ is decreasing in $[0,\ell]\subset [0,\frac{\pi}{\sqrt{2c}})$ and so
    \beq  \cos^{2}\left(\sqrt{\frac{c}{2}}\ell\right)=\alpha(0)\leq \alpha(t)\leq \alpha(\ell)=1.\eeq
By using Lemma \ref{keylinear2} and similar arguments as in the proof of Theorem \ref{main2}, we  obtain
$$
\Delta r\leq2(n-1)\frac{\mathrm{sn}_{c/2}'(r)}{\mathrm{sn}_{c/2}(r)}+\frac{\mathrm{sn}_{2c}'(r)}{\mathrm{sn}_{2c}(r)}
$$
on $M\setminus\mathrm{cut}(p)\cup \{p\}$.  If the identity in \eqref{Laplacian comparison 1} holds on $M\setminus\mathrm{cut}(p)\cup \{p\}$, then $(M,\omega_g)$ has constant holomorphic bisectional curvature $c>0$ and it is  isometrically biholomorphic to the projective space $M_c$.
\eproof

\vskip 1\baselineskip
\bproof[Proof of Theorem \ref{thm volume Ric and HSC}] Given the Laplacian Comparison Theorem \ref{thm Laplacian Ric and HSC},  the proof follows the same approach as in \cite[Theorem~1.3]{Yang25+}.
\eproof

\vskip 1\baselineskip

\section{Examples}\label{example}
We consider the model space $(\mathbb{C}\mathbb{P}^n, \omega_{\mathrm{FS}})$. It is well-known that it has
\beq \mathrm{Ric}(\omega_{\mathrm{FS}}) =(n+1)\omega_{\mathrm{FS}}, \qtq{\text{and}} \mathrm{diam}(\mathbb{C}\mathbb{P}^n, \omega_{\mathrm{FS}})=\frac{\pi}{\sqrt 2}.\eeq
 If $ r(x) = d(p,x)$ is the distance function from a fixed point $p \in \C\P^n$, then
\beq  \Delta_{\C\P^n} r= 2(n-1)\frac{\mathrm{sn}_{1/2}'( r)}{\mathrm{sn}_{1/2}( r)} + \frac{\mathrm{sn}_{2 }'( r)}{\mathrm{sn}_{2 }( r)}.\eeq
Let $(M,\omega_g)$ be a complete K\"ahler manifold with Ricci curvature $\mathrm{Ric}(\omega_g) \geq (n+1)\omega_g$. One might naturally consider comparing the diameter $d$ and the Laplacian $\Delta_M r$ of such manifolds with those of the model space. \emph{However, neither of these comparisons can hold when the complex dimension $n\geq 2$. }
 If $(M,\omega_g)$ is regarded as a Riemannian manifold, then by Myer's diameter estimate and the Laplacian comparison theorem for Riemannian manifolds, one has
\beq \mathrm{diam}(M,\omega_g)\leq \frac{\pi}{\sqrt 2} \cdot \sqrt{\frac{4n-2}{n+1}},\eeq
 and
\beq \Delta_M r\leq (2n-1)\frac{\mathrm{sn}_c'(r)}{\mathrm{sn}_c(r)}, \quad c=\frac{n+1}{2n-1}.\eeq
 When  $n\geq 2$,  these two upper bounds exceed their counterparts in the model space. Moreover, for the $n$-dimensinal  product manifold
 \beq
 (M,\omega_g) = \left( \C\P^1, \frac{2}{n+1} \omega_{\mathrm{FS}} \right) \times \cdots \times \left( \C\P^1, \frac{2}{n+1} \omega_{\mathrm{FS}} \right),
 \eeq
it has $\mathrm{Ric} (\omega_g)= (n+1)\omega_g$, and a direct calculation shows that
\beq
\mathrm{diam}(M, \omega_g) = \pi \sqrt{\frac{n}{n+1}} > \frac{\pi}{\sqrt{2}} = \mathrm{diam}(\C\P^n, \omega_{\mathrm{FS}}).
\eeq
It can also demonstrate that  the following global Laplacian comparison can not hold on $M\setminus \mathrm{cut}(p)\cup\{p\}$
\beq
\Delta_M r \leq 2(n-1)\frac{\mathrm{sn}_{1/2}'(r)}{\mathrm{sn}_{1/2}(r)} + \frac{\mathrm{sn}_{2 }'(r)}{\mathrm{sn}_{2 }(r)}.\label{model Laplacian}
\eeq
This becomes evident as  $r\nearrow \frac{\pi}{\sqrt{2}}$, where the right-hand side tends to $-\infty$. Moreover, the following example shows that the model Laplacian comparison \eqref{model Laplacian} may fail even locally. We refer to \cite{Liu11, TY12, Liu14, LY18, NZ18} for more discussions on various examples. 

\bexample \label{examplen=2}
Let $(M,\omega_g)$ be the product manifold
\beq
\left(\C\P^{1},\frac{2}{3}\om_{\mathrm{FS}}\right)\times \left(\C\P^{1},\frac{2}{3}\om_{\mathrm{FS}}\right),
\eeq
and $r$ be the distance function from a given point $p\in M$. Then for $0<r<\frac{\pi}{\sqrt{2}}$,
\beq \max_{\{x|d(x,p)=r\}} \Delta_M r(x) >2 \frac{\mathrm{sn}_{1/2}'(r)}{\mathrm{sn}_{1/2}(r)}+\frac{\mathrm{sn}_{2}'(r)}{\mathrm{sn}_{2}(r)}.\label{laplacianrelation1} \eeq
\eexample

\noindent
Before proving \eqref{laplacianrelation1}, we require the following lemma, which computes the Laplacian of the distance function on such product manifolds.

\blemma\label{lem product manifolds} Let $(M,\om_g)$ be the $n$-dimensional product manifold
\beq
\left( \C\P^1, \frac{2}{n+1} \omega_{\mathrm{FS}} \right) \times \cdots \times \left( \C\P^1, \frac{2}{n+1} \omega_{\mathrm{FS}} \right).
\eeq
 For each $1\leq i\leq n$, let $d_i$ denote the distance function in the $i$-th factor.
Fix a point $p=(p_1,\cdots,p_n)\in M$ and let $r$ be the distance function from $p$.
Then at a point $q=(q_1,\cdots,q_n)$ in $M\setminus\mathrm{cut}(p)\cup \{p\}$, we have
\beq
\Delta r(q)=\sum_{i=1}^n f_r((n+1) \lambda_i^2)+(n-1)f_r(0).\label{productlaplacian}
\eeq
where $\lambda_i = d_i(p_i, q_i) / r(q)$ and
\beq
f_r(x)=\frac{\mathrm{sn}_x'(r)}{\mathrm{sn}_x(r)}=
\begin{cases}
    \sqrt{x}\cot(\sqrt{x}r),&0<x<\pi^2/r^2,\\
    1/r,&x=0.
\end{cases}
\eeq
Moreover, $f_r(\bullet)$ is concave on $(0,\pi^2/r^2)$ for any $r>0$, and
\beq
\max_{\{x|d(x,p)=r\}} \Delta r(x)=nf_r\left(\frac{n+1}{n}\right)+(n-1)f_r(0).\label{productlaplacian1}
\eeq
\elemma

\bproof
Let $\gamma \colon [0,\ell] \to M$ be the unit-speed minimal geodesic joining $p$ and $q$, and suppose
\beq
\gamma = (\gamma_1, \dots, \gamma_n).
\eeq
Let $e_1(t), \cdots, e_{2n}(t)$ be orthonormal parallel fields along $\gamma$ such that
\beq
Je_{2i}=e_{2i-1}\qtq{and} e_{2i}=\frac{1}{\lambda_i}\gamma_i',
\eeq
for all $1\leq i\leq n$. Here, if $\lambda_i=0$ for some $i$, then $p_i=\gamma_i(t)=q_i$ for all $t\in [0,\ell]$. In this case, we choose $e_{2i}$ to be any unit vector in the tangent space to the $i$-th factor at $p_i$, and still set $e_{2i-1} = J e_{2i}$.
A direct calculation shows that the Jacobi fields along $\gamma$ are given by
\beq
V_{2i-1}(t)=\frac{\sin\left(\lambda_i\sqrt{(n+1)}t\right)}{\sin\left(\lambda_i\sqrt{(n+1)}\ell\right)}e_{2i-1}(t),\quad
V_{2i}(t)=\frac{t}{\ell}e_{2i}(t),
\eeq
where $1\leq i\leq n$. It follows that for every $1\leq i\leq n$,
\beq
\mathrm{Hess}\;r(e_{2i}(\ell),e_{2i}(\ell))=I_\gamma(V_{2i},V_{2i})= f_\ell((n+1)\lambda_i^2),
\eeq
\beq
\mathrm{Hess}\;r(e_{2i}(\ell),e_{2i}(\ell))=I_\gamma(V_{2i},V_{2i})-\int_0^\ell\langle V_{2i}'(t),\gamma'(t)\rangle^2 dt=(1-\lambda_i^2)\frac{1}{\ell}.
\eeq
Therefore, we conclude that
\beq
\Delta r=\sum_{i=1}^n f_r((n+1)\lambda_i^2)+(n-\sum_{i=1}^{n}\lambda_i^2)\frac{1}{\ell}=\sum_{i=1}^n f_r((n+1)\lambda_i^2)+(n-1)f_r(0).
\eeq
This completes the proof of \eqref{productlaplacian}.\\

Next, we shall prove that $f_r(\bullet)$ is concave on $(0,\pi^2/r^2)$ for any $r>0$. Its first derivative with respect to $x$ is
\be
f_r'(x)
&=&\frac{\cos(\sqrt{x}r)}{2\sqrt{x}\sin(\sqrt{x}r)}-\frac{r}{2\sin^2(\sqrt{x}r)}\\
&=&\frac{\cos(\sqrt{x}r)\sin(\sqrt{x}r)-\sqrt{x}r
}{2\sqrt{x}\sin^2(\sqrt{x}r)}<0,
\ee
and the second derivative is
\be
f_r''(x)
&=& -\frac{\cos(\sqrt{x}r)}{4\sqrt{x^3}\sin(\sqrt{x}r)}-\frac{r}{4x\sin^2(\sqrt{x}r)}+\frac{r^2\cos(\sqrt{x}r)}{2\sqrt{x}\sin^3(\sqrt{x}r)}\\
&=& \frac{2xr^2\cos(\sqrt{x}r)-\sqrt{x}r\sin(\sqrt{x}r)-\sin^2(\sqrt{x}r)\cos(\sqrt{x}r)}{4\sqrt{x^3}\sin^3(\sqrt{x}r)}.
\ee
If $x\in [ \pi^2/ 4r^2,\pi^2/r^2)$, then $\sqrt{x}r\in [\pi/2,\pi)$. It follows that
\be
&&2xr^2\cos(\sqrt{x}r)-\sqrt{x}r\sin(\sqrt{x}r)-\sin^2(\sqrt{x}r)\cos(\sqrt{x}r)\\
&<& \left(2xr^2-\sin^2(\sqrt{x}r)\right)\cos(\sqrt{x}r)\leq 0,
\ee
and so $f_r''(x)\leq 0$.
Thus, to show the concavity, it suffices to show
\beq
\phi(\theta)=2\theta^2-\theta\tan \theta-\sin^2\theta\leq 0,\label{concave condition}
\eeq
for all $\theta\in (0,\pi/2).$ Notice that for $\theta\in(0,\pi/2)$, we have
\beq
\sin \theta\geq \theta-\frac{1}{6}\theta^3,\quad
\tan \theta\geq \theta+\frac{1}{3}\theta^3.
\eeq
Therefore, for $\theta\in(0,\pi/2)$,
\beq
\phi(\theta)\leq 2\theta^2-(\theta^2+\frac{1}{3}\theta^4)-(\theta^2-\frac{1}{3}\theta^4+\frac{1}{36}\theta^6)=-\frac{1}{36}\theta^6.
\eeq
This completes the proof of the concavity of $f_r(\bullet)$ on $(0,\pi^2/r^2)$.
Finally, applying Jensen's inequality to the right-hand side of~\eqref{productlaplacian}, we conclude that the maximum is attained when $\lambda_i^2 = 1/n$ for all $1\leq i\leq n$. This establishes~\eqref{productlaplacian1}.
\eproof

\vone
\bproof[Proof of inequality \eqref{laplacianrelation1}]
By Lemma \ref{lem product manifolds}, we have
\begin{eqnarray}
g(r)
&:=& \max_{\{x|d(x,p)=r\}} \Delta_M r(x)-\left( 2 \frac{\mathrm{sn}_{1/2}'(r)}{\mathrm{sn}_{1/2}(r)}+\frac{\mathrm{sn}_{2}'(r)}{\mathrm{sn}_{2}(r)}\right)\nonumber\\
&=& 2\frac{\mathrm{sn}_{3/2}'(r)}{\mathrm{sn}_{3/2}(r)} + \frac{1}{r}
- 2 \frac{\mathrm{sn}_{1/2}'(r)}{\mathrm{sn}_{1/2}(r)}
- \frac{\mathrm{sn}_{2}'(r)}{\mathrm{sn}_{2}(r)} .
\end{eqnarray}
Recall that the Laurent expansion of $\cot r$ at 0 is given by
\beq
\cot r=\sum_{k=0}^\infty \frac{(-1)^k2^{2k}B_{2k}}{(2k)!}r^{2k-1},\label{taylor}
\eeq
where the convergence radius of this expansion is $(0,\pi)$, and $B_k$ are Bernoulli numbers which satisfy
\beq
B_0=1,\quad B_{2k-1}=0\qtq{and} (-1)^{k}B_{2k}<0
\eeq
for all $k\geq 1$. It follows that the Taylor expansion of $g(r)$ at $r=0$ is given
\beq
g(r)=\sum_{k=1}^\infty \frac{(-1)^k2^{2k}B_{2k}r^{2k-1}}{(2k)!}\left(\frac{3^k-1}{2^{k-1}}-2^k\right),\label{taylor2}
\eeq
where the convergence radius of this expansion is $[0,\pi)$. Since the coefficient
\beq
\frac{3^k-1}{2^{k-1}}-2^k
\eeq
vanishes when $k = 1$ and is positive for all $k \geq 2$, we conclude that $g(r)$ is a sum of positive terms. This completes the proof.
\eproof

\end{document}